\documentclass[11pt]{article}
\usepackage[french]{babel}
\usepackage{latexsym}
\usepackage{amsmath}
\usepackage{amsfonts}
\usepackage{amssymb}
\usepackage{mathrsfs}
\usepackage{eucal}
\usepackage{exscale}
\usepackage{makeidx}
\usepackage{makeidx}
\usepackage{graphicx}
\usepackage[T1]{fontenc}
\usepackage{amsmath,amssymb}
\usepackage[all]{xy}
\newtheorem{fait}{Fait}
\newtheorem{thm}{Th\'{e}or\`{e}me}%[subsection]
\newtheorem{df}{D\'efinition}%[subsection]
\newtheorem{cor}{Corollaire}
\newtheorem{lem}{Lemme}

\newcommand{\preuve}{\noindent {\bf Preuve.}\hspace{2mm}}

\newcommand{\qed}{{$\square$}\medbreak}

\begin{document}
\title{S\'eparabilit\'e en logique positive}

\author{Mohammed Belkasmi\\
Universit\'e de Lyon\\ Universit\'e Lyon 1\\
CNRS UMR 5208
Institut Camille Jordan\\
43 blvd du 11 novembre 1918\\
69622 Villeurbanne cedex, France
}
\maketitle

%\section{Introduction}
%$$  \xymatrix{
 %   T \ar@/^/[rrd] \ar@/_/[rdd] \ar@{.>}[rd] \\
 %   & A \ar[r] \ar[d] & B \ar[d] \\
 %   & C \ar[r] & D 
 % }$$

La logique positive a \'et\'e introduite par Ita\"\i\ Ben Yaacov et Bruno Poizat
dans \cite{begnacpoizat}. Dans \cite{pervers}, Bruno Poizat a \'etudi\'e
en particulier les propri\'et\'es des extensions et sous-structures positivement
\'el\'ementaires, et il a pos\'e la question suivante:
\\

\noindent
{\bf Question:} {\em Est-ce qu'une restriction \'el\'ementaire d'une structure s\'epar\'ee
est s\'epar\'ee?}
\\

Dans cette note nous proposons la preuve d'une r\'eponse positive \`a cette
question. La preuve du r\'esultat principal (le th\'eor\`eme \ref{princ})
fait usage des amalgamations dans certaines classes de structures.

\section{Logique positive: g\'{e}n\'{e}ralit\'{e}s}
La logique positive est une branche de la logique du premier ordre dont
la sp\'{e}cificit\'{e} est la non utilisation de la n\'{e}gation.
Le premier impact de cette caract\'{e}risation est la r\'{e}duction de l'ensemble des formules \`{a} l'ensemble des formules positives
%La logique positive est le fragment de la Logique du premier ordre qui n'utilise
%pas la n\'{e}gation.
% plus pr�cisement si nous consid�rons un langage L dont la signature  comprend  des
%symboles de relation de fonction et de constante individuelle.
%Les formules positives sont celles
 qu'on obtient \`{a} partir des formules atomiques gr\^{a}ce \`{a}
l'emploi des symboles ${\vee}$, ${\wedge}$, ${\exists}$.
Une formule positive se met donc sous la forme $\exists\bar{y} f(\bar{x}, \bar{y})$, 
o\`u $f(\bar{x}, \bar{y})$ est une formule libre positive 
c'est \`{a} dire sans n\'{e}gation.
Comme il n'y a pas de formule positive contradictoire, un symbole sp\'{e}cial $\bot$ d\'{e}notant l'antilogie
est ajout\'e au langage. Dans le
reste de cette section, nous rappellerons certaines d\'{e}finitions et notions
de la logique positive. Pour plus de d\'etails, \cite{begnacpoizat} est une
source suffisamment compl\`ete.

Comme dans la logique du premier ordre
avec n\'egation, un \'{e}nonc\'{e} est une formule sans variables libres. 
Un \'{e}nonc\'{e} est dit h-universel s'il est la n\'{e}gation d'un \'{e}nonc\'{e} positif ; il est donc
de la forme $\neg (\exists \bar{x}) f(\bar{x})$, 
soit encore $(\forall \bar{x})  \neg f(\bar{x})$, o\`{u} $f(\bar{x})$ est libre et positive.
La conjonction de deux \'{e}nonc\'{e}s h-universels est \'{e}quivalente \`{a} un
\'{e}nonc\'{e} h-universel ; en effet  $\neg(\exists \bar{x}) f(\bar{x}) \wedge \neg(\exists \bar{y}) g(\bar{y})$ se met sous la forme
$\neg(\exists \bar{x},\bar{y})( f(\bar{x})\vee g(\bar{y}))$. De m\^{e}me, la disjonction de deux \'{e}nonc\'{e}s h-universels est
h-universelle.

Un \'{e}nonc\'{e} positif est h-inductif simple
%; plus pr\'{e}cis\'{e}ment, 
s'il s'\'{e}crit sous la forme
\[(\forall \bar{x}) [ (\exists \bar{y})
f(\bar{x},\bar{y}) \longrightarrow (\exists \bar{z}) g(\bar{x},\bar{z}) ]\ ,\] 
o\`{u} $f$ et $g$ sont libres et positives. Sous forme pr\'{e}nexe il
devient $$(\forall\bar{u})(\exists\bar{v})(\neg\varphi(\bar{u})\vee\psi(\bar{u}, \bar{v}))\ ,$$
o\`{u} $\varphi$ et $\psi$ sont libres et positives.
% on voit
%donc que
Par cons\'equent, la disjonction de deux \'{e}nonc\'{e}s h-inductifs simples l'est aussi.
Un \'{e}nonc\'{e} h-inductif est la conjonction d'un nombre fini d'\'{e}nonc\'{e}s h-inductifs
simples ; ces \'{e}nonc\'{e}s forment une famille close par disjonction et
conjonction.

Soient $ M$ et $ N$ deux structures avec la m\^eme
signature correspondant \`a un langage L.
%Suivant
%\cite{universpositifs}, nous dirons qu'
Une application $h$ de $
M$ vers $ N$ est un homomorphisme si pour tout uple
$\bar{m}$ extrait de $M$, toute L-formule
positive $\phi$, $M\models\phi(\bar m)$ implique $
N\models\phi(h(\bar m))$. Dans ce cas, la structure $ N$ est dite
une continuation de $ M$. L'homomorphisme $h$ est dit une
immersion si $\bar m$ et $h(\bar m)$ satisfont exactement les m\^emes formules
positives dans $ M$ et $N$ respectivement.

Dans la logique positive le th\'{e}or\`{e}me de
compacit\'{e} est aussi vrai pour les \'{e}nonc\'{e}s 
qui y sont propres, et nous nous r\'ef\'ererons au th\'eor\`eme suivant comme ``compacit\'e positive''.

\begin{fait}[{\cite[Corollaire 4]{begnacpoizat}}]\label{compacite}
Une th\'eorie $h$-inductive $($ensemble d'\'{e}nonc\'{e}s h-inductifs$)$
est consistante pourvu que chacun de ses
sous-ensembles finis le soit.
\end{fait}

Deux th\'{e}ories h-inductives sont dites compagnes si elles
ont les
m\^{e}mes cons\'{e}quences h-universelles \cite{begnacpoizat}.
Le compagnonage des mod\`eles se caract\'erise
en utilisant une notion cl\'e \`a la logique positive,
d\'efinie dans \cite{begnacpoizat},
celle de mod\`ele existentiellement clos.
Un \'{e}l\'{e}ment $M$ d'une
classe $C$ de $L$-structures $C$ est dit
existentiellement clos dans $C$ si tout homomorphisme de
$M$ dans un \'{e}l\'{e}ment de $C$ est une immersion.

%Dans \cite{begnacpoizat}, un membre $M$ d'une classe de $L$-structures
%est dit positivement existentiellement clos par rapport
%\`a cette classe si
%toute continuation de $M$ dans un autre
%membre de la m\^eme classe est une immersion.

\begin{fait}[{\cite[lemme 7]{begnacpoizat}}]\label{comp}
 Deux th\'{e}ories h-inductives sont compagnes si et seulement si elles ont
les m\^{e}mes mod\`{e}les existentiellement clos.
\end{fait}

%Une cons\'{e}quence 
L'\'{e}tude dans \cite{begnacpoizat} des th\'{e}ories h-inductives
ainsi que le fait \ref{ind} ci-dessous montrent qu'une th\'{e}orie h-inductive $Ti$ a une
compagne maximale $Tk$ , qui est la th\'{e}orie h-inductive de ses mod\`eles
existentiellement clos ; les compagnes de $Ti$ sont les th\'{e}ories inductives
comprises entre sa compagne minimale $Tu$ form\'{e}e
de ses cons\'{e}quences h-universelles et sa compagne maximale $Tk$ .
La th\'eorie $Tk$ est appel\'ee enveloppe de Kaiser de Ti.
En pr\'esence des param\`etres provenant d'un
ensemble $A$, on notera $Tu(A)$ et $Tk(A)$.

\begin{fait}[{\cite{begnacpoizat}}]\label{ind}
 La classe des mod\`{e}les d'une th\'{e}orie $h-$inductive est inductive.
\end{fait}

\section{Extensions \'{e}l\'{e}mentaires en logique positive}

%\begin{df}[\cite{universpositifs}]
%Continuation/immersion aux sens d\'efinis par Poizat.
%\end{df}
La notion d'extension \'{e}l\'{e}mentaire en logique positive 
a \'{e}t\'{e} introduite et \'{e}tudi\'{e}e dans \cite{universpositifs}:

\begin{df}[{\cite{universpositifs}}]
Une continuation $N$ de $M$ est une extension \'{e}l\'{e}mentaire
de $M$, not\'ee $M\preceq_+N$,
si $N$ est un mod\`ele existentiellement clos de la classe
des mod\`eles de la th\'{e}orie $h$-universelle
$Tu(M)$ de $M$ dans le langage $L(M)$. 
\end{df}
Nous utiliserons l'expression
extension \'el\'ementaire au sens positif.

En utilisant le th\'eor\`eme 1 de
\cite{begnacpoizat},
on d\'eduit que, si $M$ est immerg\'{e} dans $N$,
ce dernier se continue en une extension \'{e}l\'{e}mentaire de $M$.

%Deux r\'esultats de Poizat et de Ben Yaacov-Poizat seront importants
%pour ce qui suit.

Dans \cite{universpositifs}, Poizat donne la caract\'erisation
suivante des extensions \'el\'ementaires en logique positive.

\begin{fait}[{\cite[Lemme 1]{universpositifs}}]\label{extensionelem}
Une continuation $N$ de $M$ en est une extension \'{e}l\'{e}mentaire
si et seulement si les deux conditions suivantes sont satisfaites:\\
1. $M$ est immerg\'{e} dans $N$.\\
2. pour tout $\bar b$ de $N$ et toute L-formule existentielle positive
$f(\bar x)$ non satisfaite par $\bar b$ dans $N$, il existe
une formule existentielle positive $g(\bar x, \bar a)$, \`{a} param\`{e}tres $\bar a$
dans $M$, qui est satisfaite par $\bar b$, et qui est contradictoire avec $f(\bar x)$ :
l'\'{e}nonc\'{e} $\neg(\exists \bar x,\bar y,\bar z)(f(\bar x,\bar y)\wedge g(\bar x,\bar z,\bar a))$,
fait partie de la th\'{e}orie universelle  $Tu(M)$,
avec $f(\bar x,\bar y)$ et $g(\bar x,\bar z,\bar a)$ des formules libres.
\end{fait}

%La preuve de ce r\'esultat d\'epend de la propri\'et\'e
%suivante qui sera cruciale dans la preuve du lemme \ref{lem2}.
%
%En fait, d'apr\`{e}s le th\'{e}or\`{e}me 23 de \cite{begnacpoizat},
%cette propri\'{e}t\'{e} caract\'{e}rise les mod\`{e}les d'une th\'{e}orie h-inductive.
%
%Une classe de structures est dite inductive si elle est close par
%limite inductive d'homomorphismes.

\section{Amalgamations}
La preuve du th\'eor\`eme \ref{princ}
utilise les liens entre
la compacit\'e et les diverses formes d'amalgamation.
Cette section est consacr\'ee \`a leur \'etude.
%Un point qui m\'erite d'\^etre soulign\'e est le lemme \ref{lem2}
%qui clarifie la question d'amalgamation des immersions.

Une des formes d'amalgamation est l'amalgamation dite asym\'{e}trique
d\'emontr\'ee dans le lemme 8 de \cite{begnacpoizat}. Le lemme
suivant en donne une version l\'eg\`erement modifi\'ee qui
convient mieux pour la suite:
%\end{fait}
 %Le lemme
%suivant en donne une version l\'eg\`erement modifi\'ee qui
%convient mieux pour la suite:

%Une autre forme d'amalgamation qui est un cas particulier de l'amalgamation asym�trique est la suivante.
\begin{lem}[{cf. \cite[Lemme 8]{begnacpoizat}}]\label{amalasym}
Soient $A$, $B$, $C$ des $L$-structures.
Si g est une immersion de $A$ dans $B$ et $h$ un homomorphisme de
$A$ dans $C$ , on peut trouver une structure $D$, mod\`{e}le de $Tk(C)$,
 un homomorphisme $g^\prime$ de $B$ dans $D$ et une immersion $h^\prime$ de
$C$ dans $D$ tels que $g^\prime\circ g = h^\prime \circ h$ .
\end{lem}
\preuve

Nommons les \'el\'ements de $A$ dans $B$ et
$C$ par les m\^emes symboles de constantes.
La preuve consiste \`a montrer que l'ensemble d'\'enonc\'es
\begin{equation}
    Tk(C) \cup diag^+(B)
\end{equation}
est consistant. En effet, si $f(\bar{a},\bar{b})$ est dans le diagramme
positif de $B$ avec $\bar{a}$ et $\bar{b}$
extraits de $A$ et de $B$ respectivement,
alors $A\models\exists \bar{y}f(\bar{a},\bar{y})$ puisque
$A$ s'immerge dans $B$. Ainsi,
on peut
interpr\'eter $\bar{b}$ par un \'el\'ement de $A$.
La formule obtenue appartient \`a $T_k(C)$.
\qed

Nous d\'emontrerons maintenant un lemme d'amalgamation
d'immersions dont le corollaire (le corollaire \ref{prince})
sera cl\'e dans la preuve du th\'eor\`eme \ref{princ}.
La forme dans laquelle ce lemme
est donn\'e ci-dessous corrige une version pr\'ec\'edente
et erron\'ee.

\begin{lem}\label{lem2}
Soient $A$ une L-structure, $B$ un mod\`{e}le de $Tk(A)$
et $C$ une $L$-structure dans laquelle $A$ s'immerge.
Alors il existe $D$ un mod\`{e}le de $Tk(C)$, et deux immersions 
$\varphi$, $\psi$, telles que le diagramme suivant est commutatif
 
 $$\xymatrix{
    A \ar[r]^{im} \ar[d]_{im} & {B} \ar[d]^{\varphi} \\
    C \ar[r]_{\psi} & {D}
  }$$
\end{lem}
\preuve
On nomme les \'{e}l\'{e}ments de $B$ et $C$ en utilisant
les m\^emes symboles pour ceux de $A$. On note $L^\star$ le langage ainsi \'elargi.
La preuve du th\'{e}or\`{e}me revient \`{a} montrer la consistance 
de la famille des  $L^\star$-\'{e}nonc\'{e}s h-inductifs

$$\Gamma=Tk(C)\cup Tu(B)\cup diag^+(B).$$ 
Soit $F=\{\chi,f(\bar\beta, \bar{b}), \neg\exists\bar y g(\bar y, \bar b) \}$  
un fragment fini de $\Gamma$, avec $\chi\in Tk(C)$, $ f(\bar\beta, \bar{b})\in diag^+(B)$
et $\neg\exists\bar y g(\bar y, \bar b)\in Tu(B)$.

 Comme $B\models \neg\exists\bar y g(\bar y, \bar b)$, et $B\models \exists \bar{x} f(\bar x, \bar{b})$,
on deduit que l'\'{e}nonc\'{e} h-inductif

$$\forall \bar z \exists \bar{x} f(\bar x, \bar{z})\rightarrow\exists\bar y g(\bar y, \bar z)$$
n'appartient pas \`a $Tk(B)$,
donc non plus \`a $Tk(A)$. Ceci implique qu'on peut trouver $\bar a\in A$ tel que
$A\models \neg\exists\bar y g(\bar y, \bar a)$, et $A\models \exists \bar{x} f(\bar x, \bar{a})$,
d'o\`u la possibilit\'{e} d'interpr\'eter nos deux \'{e}nonc\'{e}s dans $A$, par suite dans $C$,
et donc la consistance de $F$.
\qed

\begin{cor}\label{prince}
Soient $A$ une L-structure, $C$ une extension \'{e}l\'{e}mentaire positive de $A$,
et $B$ un mod\`{e}le de $Tk(A)$, alors $B$ s'immerge en un mod\`{e}le de $Tk(C)$ 
\end{cor}

\section{Restrictions \'el\'ementaires des structures s\'epar\'ees}

Dans cette section, nous d\'emontrerons le th\'eor\`eme principal
de cet article. Nous commen\c cons par le rappel des notions
principales de la preuve.

Les propri\'et\'es des types en logique positive
ont \'et\'e \'etudi\'ees dans \cite{begnacpoizat} et
\cite{pervers}.

\begin{df}[{ \cite{begnacpoizat}, \cite{pervers}}]
Soit $Ti$ une th\'{e}orie h-inductive dans un langage $L$.
 Un $n$-type est un ensemble maximal de formules
positives en $n$ variables,
qui est consistant avec $Ti$.

%De m\^{e}me
Un $n$-type \`{a} param\`{e}tres dans $M$  est un ensemble maximal de formules
positives en $n$ variables, \`{a} param\`{e}tres dans $M$,
qui est consistant avec $Ti(M)$ (ou de fa\c{c}on \'{e}quivalente, avec $Tk(M)$.)
\end{df}
On note $S_n(A)$ l'espace des n-types \`{a} param\`{e}tres dans $A$.
Un n-type de $S_n(M)$ a une r\'{e}alisation dans
une extension \'{e}l\'{e}mentaire de $M$. On d\'{e}finit
sur $S_n(A)$ une topologie dont la base
des ferm\'{e}s est l'ensemble des $F_f$ ou $f$ parcourt
l'ensemble des formules positives et $$F_f=\lbrace p\in S_n(A) \vert p\vdash f \rbrace$$

L'espace des types (positifs) est quasi-compact d'apr\`es
le fait \ref{compacite}, mais pas n\'ecessairement
s\'epar\'e. Dans {\cite{pervers}}, Poizat \'etudie des cons\'equences
du manque de s\'eparation et introduit la d\'efinition suivante:
\begin{df}
Une structure $M$ est dite s\'epar\'ee si,
pour chaque entier $n$, l'espace de types $S_n(M)$ est un quasi-compact
s\'epar\'e.
\end{df}

Dans \cite{begnacpoizat}, est d\'emontr\'ee la caract\'erisation
suivante de la s\'eparabilit\'e:
\begin{fait}[{\cite[Th\'eor\`eme 20]{begnacpoizat}}]\label{sepamal}
Les espaces de types d'une th\'{e}orie h-inductive $Ti$ sont
 s\'{e}par\'{e}s si et seulement si on peut amalgamer les
homomorphismes entre les mod\`{e}les de son enveloppe de Kaiser $Tk$ .
\end{fait}

L'\'enonc\'e suivant est aussi v\'erifi\'e dans {\cite{pervers}}:
\begin{fait}\label{extensionseparee}
Une extension \'{e}l\'{e}mentaire $N$ d'une structure s\'epar\'ee $M$ est s\'epar\'ee.
\end{fait}
Par souci de compl\'etude, nous reprenons ici sa preuve
simple qui utilise le fait suivant:
\begin{fait}[{\cite[lemme 3]{begnacpoizat}}]\label{ex}
Soit $Ti$ un ensemble d'\'{e}nonc\'{e}s h-inductifs, et soit $Tu$ l'ensemble des
\'{e}nonc\'{e}s h-universels qui sont cons\'{e}quences d'un fragment fini de $Ti$ . Alors
tout mod\`{e}le existentiellement clos de $Tu$ est mod\`{e}le de $Ti$.
\end{fait}
\preuve
Soient $N_1$, $N_2$, $N_3$, trois mod\`{e}les de $Tk(N)$ tels que $N_1$
 se continue respectivement dans $N_2$ et
dans $N_3$. Comme $N$ est une extension
\'el\'ementaire de $M$, d'apr\`es le fait \ref{ex}, $N$ est aussi
mod\`ele de $Tk(M)$.
Par cons\'equent,
$N_1$, $N_2$ et $N_3$ sont aussi mod\`{e}les de 
$Tk(M)$. En utilisant le fait \ref{sepamal} on les amalgame
dans la classe des mod\`eles de $Tk(M)$. On continue ensuite l'amalgame
ainsi obtenu dans un mod\`ele de $Tk(N)$ en utilisant
le lemme d'amalgamation asym\'etrique.
\qed
%$$\xymatrix{
% {M}\ar[r]^{\varphi_2} \ar[d]_{\varphi_3} & {N_{2}} \ar[d]^{\psi_2} \\
% {N_{3}} \ar[r]_{\psi_3} & {N'}$$

%En effet, d'apr\`es le th�or�me 20 de  {\cite{begnacpoizat}},
%les espaces de types d'une th�orie $h$-inductive $Ti$ sont
%s�par�s si et seulement si on peut amalgamer les
%homomorphismes entre les mod�les de son enveloppe de Kaiser $Tk$.
%On remarque alors que
%si on peut amalgamer dans $Tk(M)$, on peut aussi amalgamer dans la classe plus restreinte $Tk(N)$.

Le th\'eor\`eme suivant est le r\'esultat principal de ce travail
que nous \'enoncerons comme une \'equivalence
gr\^ace au fait \ref{extensionseparee}.

\begin{thm}\label{princ}
Soient $M$ et $N$ deux $L$-structures telles que $M\preceq_+N$.
Alors $M$ est s\'epar\'ee si et seulement si $N$ est s\'epar\'ee.
%Une restriction \'{e}l\'{e}mentaire d'une structure s\'epar\'ee est s\'epar\'ee.
\end{thm}

%'homomorphisme $i'$ est une immersion, donc $N_1$ est un mod\'{e}le de $Tu(N)$.\qed\
%\begin{Large}Remarque
%    \end{Large}
%Nous pouvons affblire l'hypoth\`{e}se  
%$M\preceq_+ N$ du lemme pr\`{e}c\'{e}dent on
%supposons que $N$ est mod\`{e}le de $Tk(M)$. et avoire le m\^{e}me r\'{e}sultat.

%\begin{lem}\label{lem3}
%Soient $M_1$, $M_2$ deux mod\`{e}les de $Tk(M)$
%et $\varphi$  un homomorphisme de $M_1$ dans  $M_2$,
%alors il existe $N_1$ et $M'$,  mod\`{e}les  de $Tk(N)$ et
%un homomorphisme  $\varphi '$ de $N_1$ dans  $M'$, tels que le diagramme suivant commute.
%\end{lem}
%$$\xymatrix{
%    {M_1} \ar[r]^{\varphi} \ar[d]_{i} & {M_2} \ar[d]^{i'} \\
%    {N_1} \ar[r]_{\varphi'} & {M'}
%  }$$
%\preuve

%Par hypoth\`ese, $M_1$ se continue dans $M_2$. D'apr\`es le lemme \ref{lem2},
%il existe $N_1\models Tk(N)$ tel que
%$M_1$ s'immerge dans $N_1$. Alors d'apr\`es le lemme d'amalgamation asym\'{e}trique,
%il existe une L-structure $M'$, mod\`{e}le de $Tk(M_2)$ telle que $M_2$ s'immerge dans $M'$ et
%que $N_1$ se continue dans $M'$ par le morphisme $\varphi '$.
%Le lemme \ref{lem2}
%montre qu'on peut continuer (voire immerger)
%$M'$ dans un mod\`{e}le de $Tk(N)$
%car $M'\models Tk(M)$. La commutativit\'{e} du
%diagramme r\'{e}sulte du lemme d'amalgamation.\qed

\preuve D'apr\`es le fait \ref{extensionseparee}, il suffit de
montrer la s\'eparabilit\'e de $M$ en supposant celle de $N$.
Le point principal de la preuve du th\'{e}or\`{e}me \ref{princ}
est de se ramener de $Tk(M)$ \`{a} $Tk(N)$ pour
pouvoir utiliser la propri\'et\'{e} d'amalgamation, le fait \ref{sepamal}.

Soient $M_1$, $M_2$, $M_3$ trois mod\`{e}les de
$Tk(M)$, $\varphi _2$ (resp. $\varphi _3$) homomorphisme de $M_1$ dans $M_2$
(resp. de $M_1$ dans $M_3$).
Par le lemme \ref{lem2}, il existe $N_1$,
mod\`{e}le de $Tk(N)$, tel que $M_1$ s'immerge dans $N_1$
et que le diagramme suivant commute
%Ensuite en appliquant le lemme d'amalgamation asym\'{e}trique , on obtient le diagramme commutatif suivant

%\begin{center}
%\includegraphics[scale=1]{diagramme3.png}

%\end{center}
$$\xymatrix{
    {M} \ar[r]^{im} \ar[d]_{im} & {M_1} \ar[d]^{i_1} \\
    {N} \ar[r]_{im} & {N_1}
  }$$
Comme $M_1$ s'immerge dans $N_1$, en
appliquant l'amalgamation asym\'{e}trique,
on obtient le diagramme commutatif
$$\xymatrix{
    {M_1} \ar[r]^{\varphi_2} \ar[d]_{i_1} & {M_2} \ar[d]^{i_2} \\
    {N_1} \ar[r]_{{\varphi'}_2} & {{M'}}
  }$$
avec $M'$ un mod\`{e}le de $Tk(M_2)$, et par cons\'equent
un mod\`{e}le de $Tk(M)$.
On sch\'{e}matise la construction faite jusque l\`{a} par le diagramme commutatif suivant
$$\xymatrix{
    {M} \ar[r]^{im} \ar[d]_{im} & {M_1} \ar[d]^{i_1}\ar[r]^{\varphi_2}& M_2\ar[d]^{i_2} \\
    {N} \ar[r]_{i} & N_1\ar[r]_{{\varphi'}_2}&M'
  }$$
Sur ce diagramme, on remarque que l'application
${\varphi'}_2\circ i$ d\'{e}finie de $N$ dans $M'$ est une immersion car $N$ est un mod\`{e}le
existentiellement clos de $Tk(M)$ et $M'\models Tk(M)$, ce qui implique que 
$M'$ est un mod\`{e}le de $Tu(N)$. Ceci nous permet de  continuer $M'$
dans $N_2$ un mod\`{e}le existentiellement clos de $Tu(N)$ (th\'{e}or\`{e}me 1 \cite{begnacpoizat}),  qui est aussi mod\`{e}le de $Tk(N)$ 
d'apr\`{e}s le fait \ref{ex}.
On obtient alors le diagramme commutatif suivant
$$\xymatrix{
    {M} \ar[r]^{im} \ar[d]_{im} & {M_1} \ar[d]^{i_1}\ar[r]^{\varphi_2}& M_2\ar[d]^{i_2}\ar[rd]^{f_2\circ i_2}& \\
    {N} \ar[r]_{i} & N_1\ar[r]_{{\varphi'}_2}&M'\ar[r]_{f_2}&N_2
  }$$
On refait la m\^{e}me
contruction pour $M_3$. On obtient le diagramme commutatif suivant
$$\xymatrix{
    {M} \ar[r]^{im} \ar[d]_{im} & {M_1} \ar[d]^{i_1}\ar[r]^{\varphi_3}& M_3\ar[d]^{i_3}\ar[rd]^{f_3\circ i_3}& \\
    {N} \ar[r]_{i} & N_1\ar[r]_{{\varphi'}_3}&{M''}\ar[r]_{f_3}&N_3
  }$$
o\`u $M''$ est un mod\`{e}le de $Tk(M)$, $N_3$ un mod\`{e}le de $Tk(N)$ et $f_3$ un homomorphisme. On a le diagramme suivant
$$\xymatrix{
    &M_2\ar[r]^{f_2\circ i_2}&N_2\\
   {M_1} \ar[ru]^{\varphi_{2}} \ar[rd]_{{\varphi_{3}}}\ar[r]^{i_1} & {N_1} \ar[ru]_{f_2\circ \varphi'_2}\ar[rd]^{f_3\circ \varphi'_3} \\
    &{M_3} \ar[r]_{f_3\circ i_3} & {N_3}
  }$$
Ensuite, par amalgamation des mod\`{e}les de $Tk(N)$ (le fait \ref{sepamal})
on obtient le diagramme commutatif suivant

$$\xymatrix{
    {N_1} \ar[r]^{f_2\circ\varphi^{'}_{2}} \ar[d]_{f_3\circ{\varphi^{'}_{3}}} & {{N}_2} \ar[d]^{\psi_2} \\
    {N_3} \ar[r]_{\psi_3} & {N'}
  }$$
o\`u  $N'$ est un  mod\`{e}le de $Tk(N)$, donc aussi de $Tk(M)$ (le fait \ref{ex}).
% et 
%$$\psi_2 \circ{\varphi}^{'}_{2} = \psi_3 \circ{\varphi}^{'}_{3}\ .$$
Il s'ensuit de cela que
$$\psi_2 \circ f_2\circ\varphi_{2}'
\circ i_1 =\psi_3\circ f_3 \circ\varphi_{3}'\circ i_1\ . $$
Ceci implique 
$$\psi_2\circ f_2\circ i_2\circ \varphi_2 =\psi_3\circ f_3\circ i_3\circ \varphi_3\, $$
On sch\'{e}matise cette constructuion par le diagramme commutatif suivant:
$$\xymatrix{
    &M_2\ar[r]^{f_2\circ i_2}&N_2\ar[rd]^{\psi_2}&\\
   {M_1} \ar[ru]^{\varphi_{2}} \ar[rd]_{{\varphi_{3}}}\ar[r]^{i_1} & {N_1} \ar[ru]_{f_2\circ \varphi'_2}\ar[rd]^{f_3\circ \varphi'_3} &&N'\\
    &{M_3} \ar[r]_{f_3\circ i_3} & {N_3}\ar[ru]_{\psi_3}&
  }$$
et on a le diagramme commutatif d'amalgamation suivant dans la classe
des mod\`eles de
$Tk(M)$:
$$\xymatrix{
    {M_1} \ar[r]^{\varphi_{2}} \ar[d]_{{\varphi_{3}}} & {M_2} \ar[d]^{\psi_2\circ f_2\circ i_2} \\
    {M_3} \ar[r]_{\psi_3\circ f_3\circ i_3} & {N'}
  }$$
%On sh\'{e}matise cette constructuion par le diagramme commutatif suivant:
%$$\xymatrix{
%    &M_2\ar[r]^{f_2\circ i_2}&N_2\ar[rd]^{\psi_2}&&
%    {M_1} \ar[ru]^{\varphi_{2}} \ar[rd]_{{\varphi_{3}}}\ar[r]^{i_1} & {N_1} \ar[ru]^{f_2\circ \varphi'_2}\ar[rd]^{f_3\circ \varphi'_3} &&N'\\
%    &{M_3} \ar[r]_{f_3\circ i_3} & {N_3}\ar[ru]^{\psi_3}&
%  }$$
Le th\'eor\`eme suit du fait \ref{sepamal}.
\qed
\medskip

\begin{Large}
\noindent\textbf{Remerciements}
\end{Large}

Je tiens \`{a} remercier Monsieur Altinel, Monsieur Ben Yaacov et Monsieur Poizat,
qui m'ont \'{e}clair\'{e} par leurs remarques, suggestions et corrections.
%Je tiens \`a remercier Monsieur Poizat particuli\`erement pour m'avoir
%sugg\'er\'e l'\'enonc\'e du lemme \ref{lem2}.

\bibliographystyle{plain}

\end{document}